\documentclass[sigconf]{acmart}

\renewcommand\footnotetextcopyrightpermission[1]{} 

\usepackage{multirow}
\usepackage{algorithm}
\usepackage{algpseudocode}
\usepackage{amsmath}
\usepackage{geometry}
\usepackage{algorithmicx}
\usepackage{algpseudocode}
\usepackage{color}

\usepackage{float}
\usepackage{stfloats}
\usepackage{listings}
\usepackage{xcolor}
\usepackage{tabularx}
\usepackage{threeparttable}
\usepackage{makecell}

\usepackage{pifont}
\newcommand{\cmark}{\ding{51}}%
\newcommand{\xmark}{\ding{55}}%

\AtBeginDocument{%
  }

\setcopyright{none}
\acmConference[]{}{}{}

\definecolor{mygreen}{rgb}{0,0.6,0}
\definecolor{mygray}{rgb}{0.9,0.9,0.9}
\definecolor{mymauve}{rgb}{0.58,0,0.82}
\newcommand{\cName}{\texttt}

\begin{document}

\title{StructMG: A Fast and Scalable Structured Algebraic Multigrid}


\author{Yi Zong}
\affiliation{%
  \institution{Tsinghua University}
  \city{Beijing}
  \country{China}}
\email{zong-y21@mails.tsinghua.edu.cn}

\author{Peinan Yu}
\affiliation{%
  \institution{Tsinghua University}
  \city{Beijing}
  \country{China}}

\author{Haopeng Huang}
\affiliation{%
  \institution{Tsinghua University}
  \city{Beijing}
  \country{China}}

\author{Zhengding Hu}
\affiliation{%
  \institution{University of Science and Technology of China}
  \city{Hefei}
  \country{China}}

\author{Xinliang Wang}
\affiliation{%
  \institution{Huawei Technologies Co., Ltd}
  \city{Beijing}
  \country{China}}

\author{Qin Wang}
\affiliation{%
  \institution{Huawei Technologies Co., Ltd}
  \city{Beijing}
  \country{China}}

\author{Chensong Zhang}
\affiliation{%
  \institution{Academy of Mathematics and Systems Science}
  \city{Beijing}
  \country{China}}

\author{Xiaowen Xu}
\affiliation{%
  \institution{Institute of Applied Physics and Computational Mathematics}
  \city{Beijing}
  \country{China}}

\author{Jian Sun}
\affiliation{%
  \institution{CMA Earth System Modeling and Prediction Center}
  \city{Beijing}
  \country{China}}

\author{Yongxiao Zhou}
\affiliation{%
  \institution{Tsinghua University}
  \city{Beijing}
  \country{China}}

\author{Wei Xue}
\affiliation{%
  \institution{Tsinghua University}
  \city{Beijing}
  \country{China}}
\email{xuewei@mail.tsinghua.edu.cn}





\begin{abstract}
Parallel multigrid is widely used as preconditioners in solving large-scale sparse linear systems.
However, the current multigrid library still needs more satisfactory performance for structured grid problems regarding speed and scalability.
Based on the classical 'multigrid seesaw', we derive three necessary principles for an efficient structured multigrid, which instructs our design and implementation of StructMG, a fast and scalable algebraic multigrid that constructs hierarchical grids automatically.
As a preconditioner, StructMG can achieve both low cost per iteration and good convergence when solving large-scale linear systems with iterative methods in parallel.
A stencil-based triple-matrix product via symbolic derivation and code generation is proposed for multi-dimensional Galerkin coarsening to reduce grid complexity, operator complexity, and implementation effort.
A unified parallel framework of sparse triangular solver is presented to achieve fast convergence and high parallel efficiency for smoothers, including dependence-preserving Gauss-Seidel and incomplete LU methods.
Idealized and real-world problems from radiation hydrodynamics, petroleum reservoir simulation, numerical weather prediction, and solid mechanics, are evaluated on ARM and X86 platforms to show StructMG's effectiveness. 
In comparison to \textit{hypre}'s structured and general multigrid preconditioners, StructMG achieves the fastest time-to-solutions in all cases with average speedups of 15.5x, 5.5x, 6.7x, 7.3x
over SMG, PFMG, SysPFMG, and BoomerAMG, respectively.
StructMG also significantly improves strong and weak scaling efficiencies.
\end{abstract}




\settopmatter{printacmref=false}
\maketitle

\section{Introduction}
Algebraic multigrid (AMG) is a method with optimal computational complexity in solving large-scale sparse linear system~\cite{Trottenberg}.
It is widely used as a preconditioner for iterative methods~\cite{Saad} in numerical solutions of partial differential equation (PDE) problems.
Based on multi-level grids, AMG eliminates high-frequency errors on the finer grid by smoothers, calculates and restricts residuals to the next coarser grid, and computes coarse-grid corrections recursively.
BoomerAMG~\cite{Van-Emden-HensonYang-58} from \textit{hypre} library~\cite{HYPRE-OL} is an unstructured AMG built on sparse matrices in Compressed Sparse Row (CSR) format.
However, the extra storage of integer indices and indirect memory access prevent it from achieving optimum performance for structured grid problems that are used in various applications such as weather prediction~\cite{2020QJRMS.146.3917M}, petroleum reservoir simulation~\cite{SPE_CSP_2}, and ocean modeling~\cite{struct_ocean}.
On structured grids, nonzero entries in a specific row of the matrix correspond exactly to the neighbors of an element in the grid.
Hence, \textit{hypre} provides SMG~\cite{BrownFalgout-127, Schaffer-125, FalgoutJones-126}, PFMG~\cite{FalgoutJones-126, AshbyFalgout-124}, and SysPFMG (system version of PFMG for vector PDE\footnote{Vector PDE: multiple unknowns associated with each grid element.}) that make use of structured grids and stencils to build hierarchical grids algebraically.
They are domain-specific AMG aiming for higher performance in structured grid problems.

\label{txt:StructInfo}

\begin{figure}[h]
  \centering
  \includegraphics[width=\linewidth]{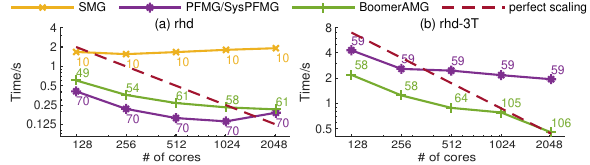}
  \caption{Strong scaling tests of \textit{hypre} 's structured AMG against unstructured AMG as preconditioners for rhd problems. Numbers of iterations are labeled near the lines. }
  \label{fig:sota}
  \vspace{-1.0em}
\end{figure}

However, none of these multigrids could perform comprehensively well enough for large-scale structured grid problems in terms of speed and scalability in various application tests. 
Take strong scaling tests of a radiation hydrodynamics (\cName{rhd}) problem discretized on a grid of $128^3$ as an example in Figure~\ref{fig:sota} (the same results with breakdowns are shown in Figure ~\ref{fig:strong}, and detailed specifications refer to Section~\ref{txt:problems}). 
For the scalar PDE\footnote{Scalar PDE: only one unknown associated with each grid element.} case in (a), SMG achieves the desired accuracy with the fewest number of iterations (denoted by \#iter in the following text) but with the longest wall time and worst scalability;
PFMG is the fastest but fails to scale continuously;
BoomerAMG scales slightly better than PFMG, but its \#iter increases with the number of cores.
Fewer \#iter of SMG and faster speed of PFMG demonstrate the advantages of structured information utilization, but there is still room for improvement in their scalability.
For the vector PDE case in (b), the results deviate from expectations that SysPFMG is significantly slower than unstructured AMG.
The current deficiency for structured problems is admitted in a technical report~\cite{osti_2021}.
A natural expectation is that a faster and more scalable structured multigrid should be well-addressed.

But improving speed and scalability at the same time remains a challenge.
Apart from the setup phase, speed depends on \#iter and solving time per iteration that is approximately proportional to the memory volume to access in each iteration. Scalability depends on suitability to utilize parallelism and variance of \#iter with increasing parallelism. 
A sophisticated method, such as plane smoothing in SMG~\cite{Schaffer-125}, generally keeps \#iter stable but has a more significant overhead for each iteration. 
Easy-to-parallelize methods, such as weighted Jacobi and Red-Black Gauss-Seidel (RB-GS) smoothers in PFMG~\cite{AshbyFalgout-124} and hybrid Jacobi/GS smoother in unstructured AMG~\cite{Van-Emden-HensonYang-58}, often have poorer convergence. 
Better approaches need to balance these conflicts.

Towards faster speed and better scalability, higher dimensional coarsening than SMG, PFMG, and SysPFMG is supported in our StructMG to reduce grid and operator complexities.
Lower complexity means less computation, less memory to access, and fewer numbers to communicate.
Robust yet efficient smoothers are implemented in a multi-thread framework to accelerate convergence, including dependence-preserving point-wise Gauss-Seidel (PGS), line-wise GS (LGS), and incomplete lower-upper (ILU) factorization.

Specifically, this paper makes the following contributions.
\begin{itemize}
  \item We design StructMG, a structured AMG that automatically constructs the hierarchy from the original matrix.
  \item We propose a symbolic analysis and code-generation method for multi-dimensional Galerkin coarsening to reduce grid and operator complexities. It enables the flexible generation of fusioned triple-matrix product with arbitrary patterns.
  \item We implement a unified sparse triangular solver (SpTRSV) framework to support robust yet efficient smoothers on structured grids. It makes a good trade-off between low cost per iteration and convergence with increasing parallelism.
  \item We evaluate idealized and real-world problems on ARM and X86 platforms. StructMG obtains the best time-to-solutions in all tests with average speedups of \textbf{15.5x}, \textbf{5.5x}, \textbf{6.7x}, and \textbf{7.3x} over \textit{hypre}'s SMG, PFMG, SysPFMG, and BoomerAMG, respectively.
  StructMG also significantly improves strong and weak scaling efficiencies.
\end{itemize}

\section{Background}\label{sec:background}
\begin{figure}[htbp]
  \centering
  \includegraphics[width=0.9\linewidth]{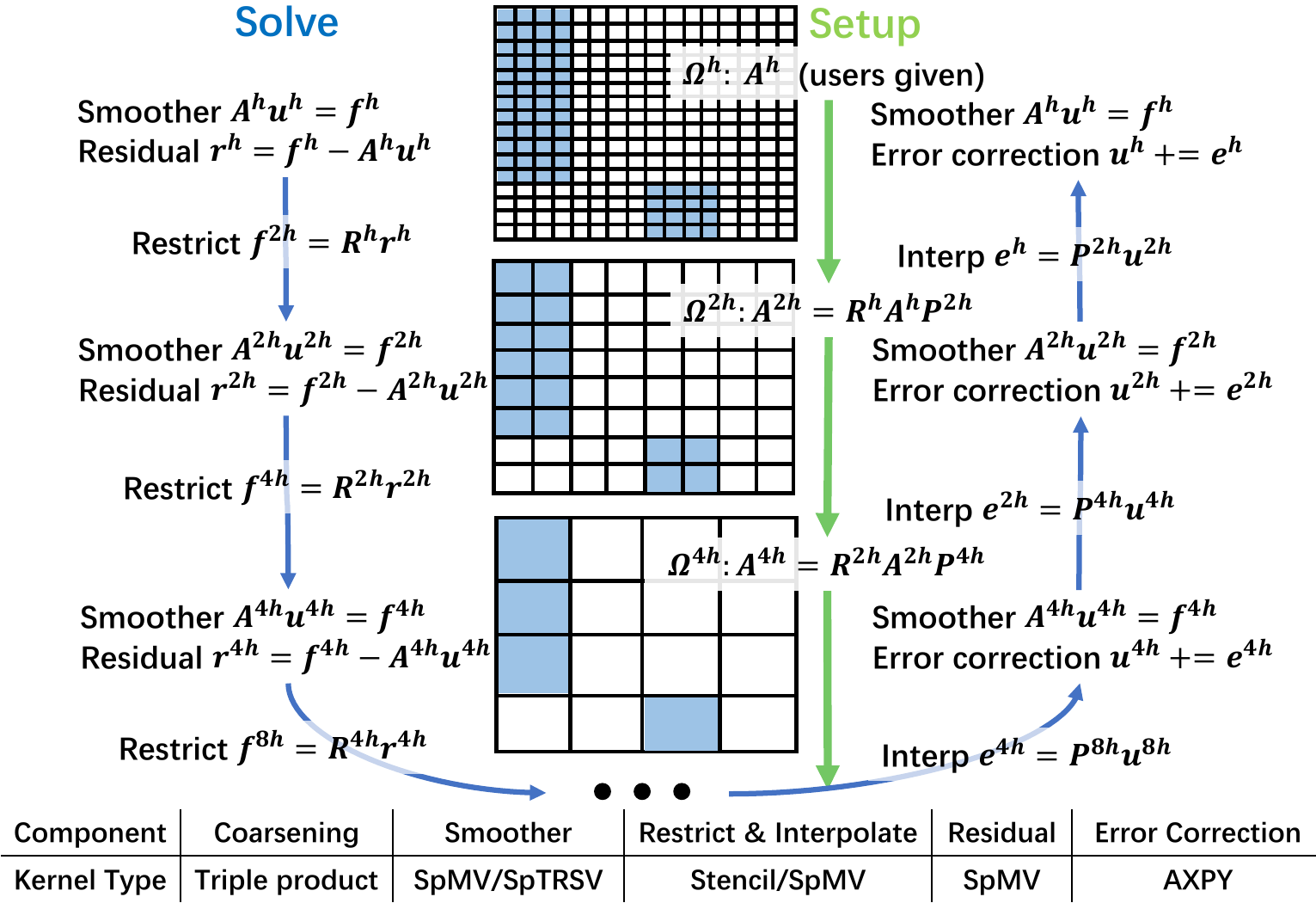}
  \caption{AMG overview. V-Cycle in the solve phase~\cite{ZY_FP16}.}
  \Description{...}
  \label{fig:ov_MG}
    \vspace{-1.0em}
\end{figure}

An illustrative overview of AMG is shown in Figure~\ref{fig:ov_MG}.

Coarsening is the critical point in the setup phase, denoted by green arrows in Figure~\ref{fig:ov_MG}.
AMG constructs a series of coarser grids $\Omega^{2h}, \Omega^{4h}, ...$, based on the finest-level grid $\Omega^h$, where the superscripts denote grid spacings of levels. 
Except for the finest level, AMG computes the coarser operator (i.e., matrix) $A^{2^{k}h}$ based on $A^{2^{k-1}h}$ of the finer level.
The automatic setup phase without users' explicit involvement is what distinguishes AMG from geometric multigrid (GMG), and also constitutes its key advantage~\cite{Stuben2000AlgebraicM}.

In the solve phase, denoted by blue arrows, AMG usually executes a V-Cycle~\cite{Trottenberg}, which starts from the finest level, traverses to the coarsest, and then reverses.
On each level, AMG invokes a smoother to solve $A^{*h}u^{*h}=f^{*h}$ approximately, and computes residual $r^{*h}$ in the downward pass;
computes error correction to update solution $u^{*h}$, and
then invokes the smoother again in the upward pass. 
Restrictions of residuals $r^{*h}$ and interpolations of errors $e^{*h}$ occur
in downward and upward pass, respectively.

The multigrid is set up once and applied iteratively. Therefore, the total time $T_\text{tot}$ of employing a multigrid in iterative solvers is 
\begin{equation} \label{eq:time}
T_\text{tot}=T_{\text{setup}} + \text{\#iter} \cdot T_{\text{single}} 
\end{equation}
where $T_{\text{setup}}$ is setup time, $T_{\text{single}}$ is single-iteration solve time, and \#iter is number of iterations to solve the linear system.
These three terms are closely related to each other.
Two primary metrics to evaluate the computational overhead of a multigrid are grid complexity $C_{\text{G}}$ and operator complexity $C_{\text{O}}$~\cite{Stuben2000AlgebraicM},
\begin{equation} \label{eq:complexity}
C_{\text{G}} = \frac{\sum_l n_l}{n_0} \quad \text{and} \quad C_{\text{O}} = \frac{\sum_l Z_l}{Z_0}
\end{equation}
where $n_l$ and $Z_l$ denote the number of unknowns and nonzero entries, respectively, on level $l$. Note that $l=0$ corresponds to the finest level given by user.
Particularly, $C_{\text{O}}$ is positively correlated with $C_{\text{G}}$ and the average number of nonzero entries per row.
The higher these two complexity metrics are, the higher the computational cost per AMG call is required.

Coarsening and smoothers are two essential components differentiating a multigrid library from others.
The other components are mostly similar, whose kernel types are shown in Figure~\ref{fig:ov_MG}.
Specifically, multigrid libraries usually apply smoothers in a block-Jacobi way across processes to enable efficient inter-process parallelism, and their differences lie in the intra-process algorithms.
The multigrids from \textit{hypre}~\cite{HYPRE-OL}, the state-of-the-art (SOTA) high-performance library of preconditioners and solvers, will be reviewed with respect to coarsening and smoothers in the next section.

\section{Related Work} \label{sec:related}

\begin{figure}[htbp]
  \centering
  \includegraphics[width=0.9\linewidth]{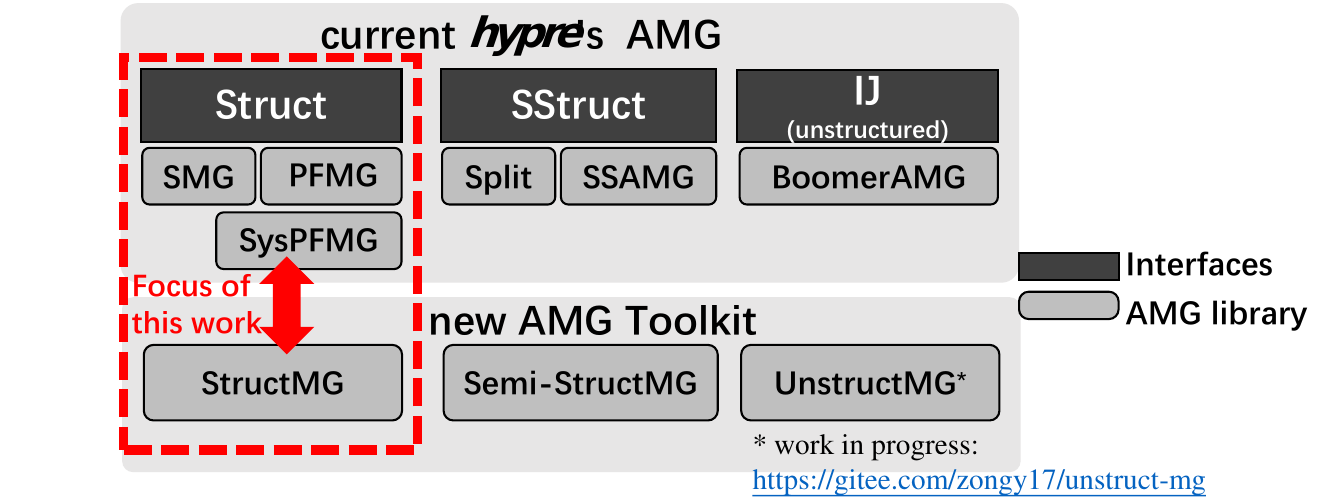}
  \caption{Illustration of the correspondence between the AMG toolkit and the different interfaces of \textit{hypre}. This article, StructMG, aims for a better solution for \textbf{Struct} scenarios.}
  \Description{...}
  \label{fig:correspond_to_hypre}
  \vspace{-1.0em}
\end{figure}

\textit{Hypre}~\cite{HYPRE-OL} enjoys a prestigious reputation in the HPC community, offering a suite of AMG libraries tailored for diverse application scenarios.
They belong to distinct interfaces, namely the Structured-Grid System (\textbf{Struct}) Interface, the Semi-Structured-Grid System (\textbf{SStruct}) Interface, and the Linear-Algebraic System (\textbf{IJ}) Interface, as shown in Figure~\ref{fig:correspond_to_hypre}.
The scenarios they cater to are suggested by the names of these three interfaces.

Our review starts with BoomerAMG, the most popular and widely-used library in \textit{hypre}.
As an unstructured implementation of algebraic multigrid under the \textbf{IJ} Interface, BoomerAMG is a \textbf{black-box} preconditioner that only relies on information provided by CSR matrices. 
BoomerAMG is particularly robust for scalar PDEs.
It selects coarse elements more flexibly with a larger $T_{\text{setup}}$ than its structured-grid-specific counterparts.
The extra overhead of indirect memory access with CSR format increases $T_{\text{single}}$, especially when 8-byte BigInt is necessary.
Its typical bottlenecks include the triple-matrix product in Galerkin coarsening in the setup phase~\cite{XuXiaowen_243}, and extensive communication requirements in the solve phase~\cite{GahvariBaker-71}. 
Many researchers focused on alleviating these hotspots~\cite{Falg_nonGal, Treister2015NonGalerkinMB, Bienz2015ReducingPC, 10.1109/IPDPSW.2013.164, 10.1177/1094342020925535, Park-SC15, Xu-aAMG}.
Comprehensive optimizations of BoomerAMG can be found in~\cite{Park-SC15}.
Its default smoother is hybrid Jacobi/GS~\cite{BoomerAMG-OL}, in which inter-thread (i.e., off-core) dependence is Jacobi-type while intra-thread (i.e., on-core) dependence is GS-type.

The other multigrids under the \textbf{Struct} and \textbf{SStruct} Interfaces are \textbf{grey-box} methods that build on structured and semi-structured grids, respectively.
A semi-structured grid consists of multiple structured parts that can be distorted, stretched, rotated, etc., and can have arbitrary connections between them.
Therefore, SSAMG~\cite{SSAMG} and Split~\cite{Split-OL} both utilize PFMG or SMG as the underlying component within each structured part.
The multigrids under the \textbf{Struct} Interfaces have an (off)diagonal storage scheme free of indirect data access \cite{osti_2021}, and provide 1D semi-coarsening for
equations discretized by finite-difference, finite-volume, and finite-element methods on logically rectangular grids.
Currently, these structured-grid-specific multigrids only support matrices of 2d5, 2d9, and 3d7, 3d19, 3d27 nonzero patterns shown in Figure~\ref{fig:nnz_layout}.

SMG~\cite{BrownFalgout-127, Schaffer-125, FalgoutJones-126} coarsens in the $z$-direction and uses a plane smoother. 
On the $xy$-plane, its smoother employs one cycle of a 2D SMG method, which in turn coarsens in $y$ and uses $x$-line smoothing. 
The heavyweight plane smoother requires both long $T_{\text{setup}}$ and $T_{\text{single}}$, but reduces \#iter remarkably.
Frequent communications in the plane smoother make SMG the poorest scalable among all multigrid methods.
\textit{Hypre}'s recent report~\cite{osti_2021} pointed out that it requires many recursions and has many short loops, which are unsuitable for parallel computations.

PFMG~\cite{FalgoutJones-126, AshbyFalgout-124} uses point-wise smoothers (like weighted Jacobi and RB-GS) instead. This is the main difference between PFMG and SMG.
As a result, PFMG is less robust but usually faster in each iteration and has better scalability.
PFMG supports both Galerkin and non-Galerkin coarsening.
The non-Galerkin approach maintains stencil patterns of 5-point in 2D problems and 7-point in 3D problems, reduces cost, and improves performance, but is usually less convergence-favored. 
As reported in~\cite{osti_2021}, PFMG is sometimes not sufficiently robust even for structured problems.
Users may turn to BoomerAMG when PFMG is too slow due to too many \#iter.
But BoomerAMG is usually not the best, as shown in Figure~\ref{fig:sota}(a).

SysPFMG\footnote{We classify SysPFMG under \textbf{Struct} Interface although it is implemented through the \textbf{SStruct} interface, because it can be used only for problems with one structured part.}
~\cite{Structured-MG-OL}
is a direct generalization of PFMG for vector PDEs.
It inherits the strengths as well as weaknesses of PFMG.
SysPFMG only supports the Galerkin approach to construct coarse-grid operators.
Interpolations are within the same variable, and nodal-type relaxations are used.
The lightweight point-wise smoothers in PFMG and SysPFMG are almost free of setup overhead and easy to parallelize.
But their \#iters may increase considerably and thus slow down the overall speed.
An aggressive option of skipping relaxation is provided in PFMG and SysPFMG but is less convergence-favored.

A comprehensive evaluation of the above \textit{hypre}'s multigrids for Laplace and Poisson equations with up to 125,000 cores can be found in~\cite{Baker2012}, and a more recent technical report~\cite{osti_2021}.
Black-box multigrid~\cite{boxMG_Factor3,boxMG_500kCores} is another structured AMG, but it has only been tested on 2D problems.
Since its source code or binary library is not public available, and our testing problems are in 3D, we do not include it in our experiments.
%

As illustrated in Figure~\ref{fig:correspond_to_hypre}, this work on StructMG corresponds to the multigrids within the \textbf{Struct} Interfaces of \textit{hypre}.
Just as PFMG is utilized by SSAMG to construct semi-structured multigrid, StructMG not only provides high-performance solutions directly for structured grid problems, but also acts as a vital building block in the success of the more complex Semi-StructMG~\cite{Semi-StructMG}.
The GPU version of this work has completed, but details and results are omitted here due to page limit, which will be in a follow-up paper.
Prior to developing a complete library on GPU, this work analyzes fundamental design principles, establishes the algorithmic framework, and validates their effectiveness.
Given the demonstrated efficiency advantages of structured kernels on GPU~\cite{Stencil–CSR_GPU}, we expect structured AMG to achieve significant performance gains on GPU.

\section{StructMG}

The related work in Section~\ref{sec:related} implies that designing a fast and scalable multigrid requires balancing different components. 
In this section, we will first introduce our design principles of StructMG, which are derived from the observations of the SOTA libraries.
Following that, two critical implementations that enable the key functionality and performance of StructMG will be discussed.

\subsection{Basic Principles of Design} \label{txt:principles}

Our design is motivated by the following quotation from the highly-influential AMG literature~\cite{Stuben2000AlgebraicM}:\\[0.5em]
\textit{\textbf{Efficient interplay between smoothing and coarse-grid correction is required for any multilevel approach.}}\\[0.5em]
We refer to this as the \textbf{multigrid seesaw}.
The different emphases on the two ends drive the distinct evolution of AMG and GMG, as analyzed in~\cite{Stuben2000AlgebraicM}.
This conclusion still holds for the contrasting design philosophies of unstructured AMG and structured AMG.
Putting most efforts on the coarse-grid correction end of the seesaw, unstructured AMG usually fixes the smoother to relatively simple relaxation schemes such as point Jacobi or GS relaxation, and enforces the efficient interplay with the coarse-grid correction by choosing the coarser elements appropriately.
Its coarser element selection is dynamic and adaptive for matrix coefficients.
However, the dynamic selection can destroy the inherent properties of structured grids, where the number of neighbors and their displacements of grid points are fixed.
Consequently, the performance penalties in unstructured scenarios, such as the extra storage of integer indices, indirect memory access, and more complicated communication pattern, have to be paid.

An effort made by PFMG to address coarse-grid correction is 1D semi-coarsening, which coarsens the elements with a fixed stride of 2 in the direction where the matrix coefficients exhibit the strongest coupling.
However, this effort is limited that it only works when the anisotropy of the matrix coefficients is grid-aligned.
Even worse, it can result in a high complexity $C_{\text{G}}=1+1/2+1/4+\cdots\rightarrow2$ for the generated multigrid.
For an evident contrast, unstructured AMG can often reduce the complexity to below 1.3~\cite{ZY_FP16} through aggressive coarsening~\cite{aggressive_coarsen}.
While PFMG shows weaker effect in coarse-grid correction than unstructured AMG, it does not invest more efforts on smoothing, the other end of the seesaw.
Its Jacobi and RB-GS smoothers have similar convergent properties with the hybrid Jacobi/GS counterparts in unstructured AMG.

Based on the above observation, to fully exploit the performance potential that structured problems offer to the multigrid method, our design of StructMG must adhere to the following principles:
\begin{description} 
    \label{principle}
    \item[P1] Data structures and implementations should be compatible with structured grids for optimum kernel performance.
    \item[P2] Multi-dimensional coarsening should be applied to reduce grid and operator complexities for faster $T_\text{setup}$ and $T_\text{single}$.
    \item[P3] Robust-yet-efficient smoothers should be used to compensate for fix-stride coarsening to obtain good convergence.
\end{description}

The first principle, \textbf{P1}, is inherited from current structured multigrids such as SMG, PFMG, and SysPFMG.
Problems on structured grids have a characteristic that there is the same number of nonzero entries in each row of the discretized matrix, and they correspond exactly to the neighbors of an element
\footnote{In this section, the term \textbf{\textit{element}} is used to refer to both vertex and cell since our method applies for both vertex-centered~\cite{KhalilWesseling-103, Dendy_BlackBox} and cell-centered~\cite{Wesseling-101, Kwak_Cell_Centered} multigrids.}
in the grid.
Therefore, the structured-grid-diagonal (SG-DIA) format~\cite{SGDIA}, which stores vectors and matrices in a multi-dimension-array way, could be exploited to avoid extra integer indices arrays.
We refer to a matrix as a structured matrix if its nonzero pattern fits the SG-DIA format.
Specifically, a structured matrix of $N$ rows and $nzpr$ entries per row with 0-valued padding at boundaries occupies exactly $O(N\times nzpr)$ memory space.
Note that in this way, a structured matrix can be formulated into a stencil with variable coefficients.
Computation and communication of SpMV and SpTRSV on structured grids will be stencil-like, as in~\cite{Wang-ICPP18}.
Avoiding the extra overhead of the general sparse matrix format is beneficial to shorter $T_{\text{single}}$, as PFMG has shown in Figure~\ref{fig:sota}(a).

Principle \textbf{P2} aims to overcome the issue of excessively high complexity in current structured multigrids.
Rather than coarsening in only one specific direction, our multi-dimensional coarsening could reduce $C_\text{G}$ to $1 + 1/8 + 1/64 + \cdots \rightarrow 8/7$ in 3D problems.
Lower grid and operator complexities reduce $T_{\text{setup}}$ and $T_{\text{single}}$ by decreasing computational workload, memory access and communication volumes.
Moreover, since all components on each level require halo updates, the lower $C_\text{G}$ results in fewer AMG levels, leading to significantly fewer communication rounds and improved scalability.
Section~\ref{txt:galerkin} will detail this indispensable technology.

Principle \textbf{P3} states that more efforts are required on the smoothing end of the seesaw.
Since fix-stride coarsening is the only possible way to maintain structured characteristic on coarser levels and must be adopted for optimum performance, its convergence weakness should be carefully compensated.
Carefulness means that one should not overlook efficiency in the excessive pursuit of convergence.
Table~\ref{tab:smoother_list} lists the smoothers in different multigrids and summarizes their convergence and efficiency characteristics.
Different from the hybrid Jacobi/GS and RB-GS that break dependence
\footnote{Strictly speaking, coloring in RB-GS is a form of reordering, but from the perspective of the original ordering, it disrupts the dependency order.
The convergence of this reordering is inferior to that of the original ordering.}
among elements, all GS-type smoothers in StructMG respect the original~dependency order and thus, are more convergence-favored.
In addition to the commonly-used point-wise smoothers, StructMG provides more advanced line-wise smoothers that are suitable for grid-aligned anisotropic problems.
For skewed or arbitrarily-oriented anisotropy that is not grid-aligned, a flexible ILU smoother is suitable to improve convergence.
The flexibility will be explained in Section~\ref{txt:StructMG_feat}.
The better convergence of these smoothers will be examined by ablation experiment in Section~\ref{txt:abla_convergence}.
All smoothers in StructMG support multi-threading to maximize the core and bandwidth utilization of modern architecture, and their efficiency will be demonstrated in Section~\ref{txt:abla_kernel}.

\begin{table}[htbp] \footnotesize
\caption{Comparison of smoothers in different multigrids. 'MT' short for multi-threading support.} \label{tab:smoother_list}
\begin{threeparttable}
  \begin{tabular}{c|cccc}
    \toprule[0.2em]
    \textbf{AMG} & \textbf{Smoother} & \textbf{MT} & \textbf{Convergence} & \textbf{Efficiency} \\
    \midrule[0.2em]
    \multirow{2}*{BoomerAMG\tnote{1}} & Jacobi, hyb. Jacobi/GS & \cmark & Ordinary & High \\
                             & ILU  & \xmark & Good & Ordinary \\
    \hline
    SMG & Plane smoothing & \cmark & Very Good & Low \\
    \hline
    PFMG/SysPFMG & Jacobi, RB-GS & \cmark & Ordinary & High \\
    \hline
    \multirow{3}*{StructMG} & Point Jacobi, GS & \cmark & Ordinary & High \\
                            & Line Jacobi, GS & \cmark & Good & High \\
                            & ILU & \cmark & Good & High \\
    \bottomrule[0.2em]
  \end{tabular}
\begin{tablenotes}    
\footnotesize               
\item[1] Other smoothers are available in BoomerAMG. Only the most commonly used ones are listed here.
\end{tablenotes}
\end{threeparttable}
\vspace{-1.0em}
\end{table}

\begin{figure}[htbp] \centering
  \includegraphics[width=0.75\linewidth]{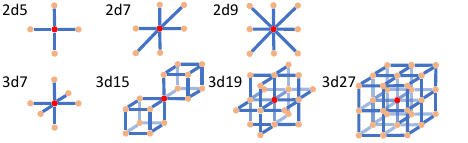}
  \caption{Typical nonzero patterns of structured matrices.}
  \label{fig:nnz_layout}
\vspace{-1.0em}
\end{figure}

\subsection{Overview of StructMG Features}\label{txt:StructMG_feat}
This subsection provides a high-level sketch of StructMG.
\subsubsection{Parallelization}
The parallelization on hierarchical grids across processes is in a domain-decomposition fashion. 
Each process takes a portion of the global grid with halo regions and collaborates with others in a bulk synchronous parallel (BSP) model.
As shown in Figure~\ref{fig:ov_MG}, the blue regions represent the multiple boxes a process may take. 
The process partitioning method is similar to \textit{hypre}.
The two critical components, coarsening and smoothers, are analyzed in the following subsections.
The other components are similar to a regular AMG implementation and will not be explained in detail.
All these components are parallelized by hybrid MPI/OpenMP.

\subsubsection{Input Data Format}
StructMG requires the same input data as SMG, PFMG, and SysPFMG. 
The input only includes the finest-level matrix $A$ and right-hand-side vector $b$, which makes a convenient and user-friendly interface possible.
The user interface of StructMG can be referred to the user guide of Semi-StructMG~\cite{Semi-StructMG}, as the former is fully integrated within the latter.
When a semi-structured grid has only 1 part, it simplifies to a structured grid.

\subsubsection{Extensibility \& Flexibility}
From its inception, StructMG was designed not to be confined to supporting structured matrices with only a few specific nonzero patterns, although the types it currently supports, as shown in Figure~\ref{fig:nnz_layout}, already cover a wide range of application scenarios.
Each type of nonzero patterns can be encoded as a mask, where 0 and 1 represent the absence or presence of a neighbor at that position, respectively.
The currently supported mask types are those with neighbor offsets in range of \{-1,0,1\} in three directions.
While keeping the multigrid framework unchanged, adding a new pattern type only requires implementing the corresponding SpMV, SpTRSV, and triple-matrix-product kernel interfaces to enable that type.
Therefore, StructMG possesses adequate extensibility to accommodate more types.
The flexibility of the ILU smoother is also reflected in the control of masks.
Users can specify the positions of nonzero entries retained by the ILU factorization through different masks, which are not necessarily the same as the mask of the original structured matrix.
For example, a 3d19 matrix can be factorized in the ILU(0) manner with 3d19 mask specified, or with more fill-ins using 3d27, or with fewer fill-ins using 3d7.

\subsection{Multi-dimensional Galerkin Coarsening: Stencil-based Triple-Matrix Product} \label{txt:galerkin}

Galerkin coarsening is commonly used in AMG to compute coarse-level matrix $A^C$.
Given a fine-level matrix $A^F$, interpolation matrix $P$, and restriction matrix $R$, Galerkin coarsening calculates the following triple-matrix product:
\begin{equation}
    A^C = RA^FP
\end{equation}
where $P, R$ are sparse structured matrices used in interpolation and restriction, respectively. Their patterns differ in different coarsening algorithms.
The challenge of multi-dimensional coarsening in structured grids lies in the intricacy of its implementation.
A naive way is to explicitly assemble $A^F$, $P$, and $R$ into CSR format and then call the general sparse matrix-matrix multiplication (SpGEMM) procedure. 
But it is not optimum for structured-grid problems since the data structure violates the principle \textbf{P1}.

The native data structure of structured grids offers the chance of faster speed and better scalability.
Hence, StructMG implements Galerkin coarsening based on grids and stencils.
Note that nonzero locations of a structured matrix (i.e., offsets of neighbors in grids) are the same across different rows (i.e., elements in grids).
This feature makes symbolic analysis and automatic generation suitable because the formulas to fuse the two-step products (i.e., which three nonzero entries could be multiplied and accumulated to the final target) could be derived, generated and compiled once, and then reused by matrices with the same patterns but different values.

\begin{figure}[htbp] \centering
  \includegraphics[width=0.9\linewidth]{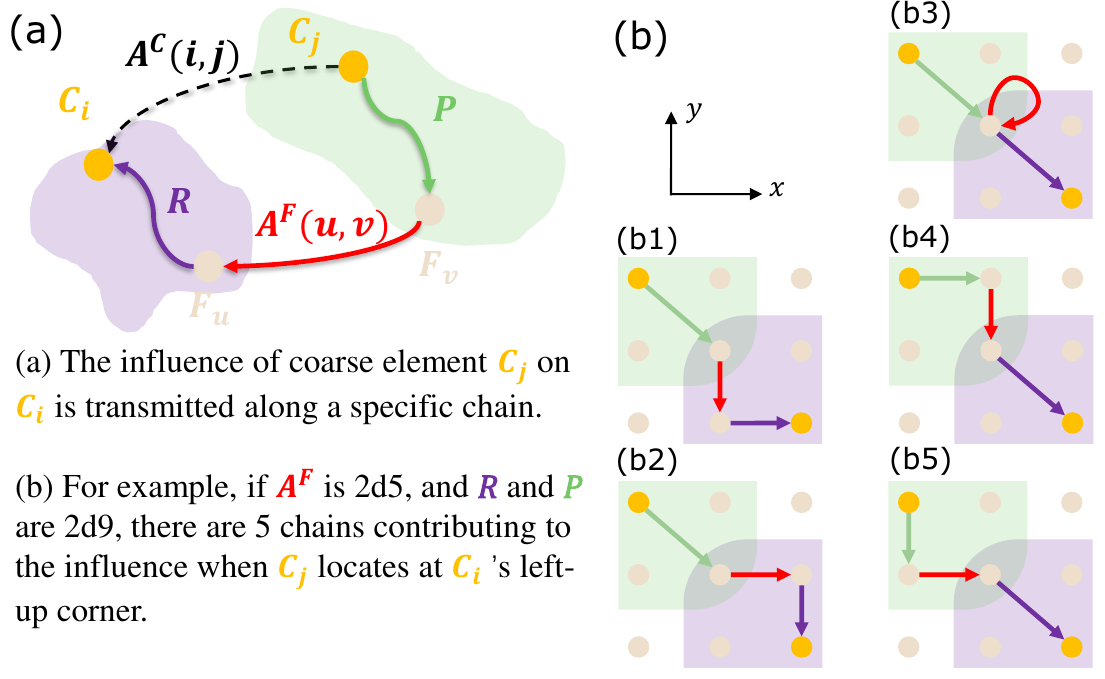}
  \caption{Illustration of stencil-based triple-matrix product $RA^FP$ on structured grids.}
  \label{fig:Galerkin-overview}
\vspace{-1.0em}
\end{figure}

To this end, we first consider the triple-matrix product from a different perspective. Since a matrix is essentially a linear operator, the product of multiple matrices can be viewed as the composition of linear operators.
With this insight, each nonzero entry in the coarse matrix $A^C$ could be disassembled into a sum of scalar products from matrix $R$, $A^F$, and $P$.
Specifically, $A^C(i,j)$, the nonzero entry at $i$-th row and $j$-th column of the coarse matrix, represents the influence of coarse element $C_j$ on coarse element $C_i$. 
Figure~\ref{fig:Galerkin-overview}(a) shows that a portion of the influence is transmitted along a specific chain through fine elements $F_v$ and $F_u$.
The green area, where $F_v$ resides, is the range of fine elements that can be influenced by the coarse element $C_j$ via interpolation.
The purple area, where $F_u$ resides, is the range of fine elements that can impose influence on coarse element $C_i$ via restriction.
$F_v$ is exactly a fine-grid neighbor influencing $F_u$ with the value of $A^F(u,v)$. Thus, as shown in the following formula, end-to-end influence is accumulated to $A^C(i,j)$. 
\begin{equation} \label{eq:RAP} \small
  A^C(i,j) = \sum_{u,v} R(C_i,F_u)A^F(u,v)P(F_v,C_j)
\end{equation}
The above equation changes the perspective from computing matrix product to finding influence pathways, but only improves performance when structured information is utilized.
As mentioned previously, $A^F$, $P$, and $R$ can be expressed as stencils with variable coefficients on structured grids.
Once their stencil patterns are given, the pattern of coarse matrix $A^C$ is determined as analyzed in~\cite{MohrWienands-104}.
The stencil coefficients at each coarse element, namely nonzero values of $A^C$ in each row, can use a common set of formulas derived statically to compute.
For the example in Figure~\ref{fig:Galerkin-overview}(b), the formulas in (b3), (b4), and (b5) are derived as follows.
Given the coordinates $C_i=\left[X,Y\right]^T$ and $C_j=\left[X-1,Y+1\right]^T$, the range of purple area, denoted by $R_{\text{from}}$, in the fine grid is 
\begin{equation} \small
R_{\text{from}}  = 
\begin{bmatrix} (X-h_x)s_x + b_x + h_x \\ (Y-h_y)s_y + b_y + h_y\end{bmatrix}  + \text{Offset(2d9)}
\end{equation}
and the range of green area, denoted by $P_{\text{to}}$, is
\begin{equation} \small
P_{\text{to}} = 
\begin{bmatrix} (X-1-h_x)s_x + b_x + h_x \\ (Y+1-h_y)s_y + b_y + h_y \end{bmatrix} + \text{Offset(2d9)}
\end{equation}
where $s_x$, $s_y$ are coarsening strides, $h_x,h_y$ are halo widths, and $b_x,b_y$ are base offsets of the coarse grid relative to the fine grid in 2D, and 2d9 offsets are
\begin{equation}\small
\text{Offset(2d9)} = \left\{
\begin{bmatrix} -1 \\-1 \end{bmatrix},
\begin{bmatrix} -1 \\ 0 \end{bmatrix},
\begin{bmatrix} -1 \\ 1 \end{bmatrix},
\begin{bmatrix}  0 \\-1 \end{bmatrix},
\begin{bmatrix}  0 \\ 0 \end{bmatrix},
\begin{bmatrix}  0 \\ 1 \end{bmatrix},
\begin{bmatrix}  1 \\-1 \end{bmatrix},
\begin{bmatrix}  1 \\ 0 \end{bmatrix},
\begin{bmatrix}  1 \\ 1 \end{bmatrix} \right\}
\label{eq:off_2d9}
\end{equation}
For each fine element $F_u \in R_{\text{from}}$, neighbors of $F_u$ are found by 2d5 offsets of fine matrix as
\begin{equation}\footnotesize
    \text{Ngb}(F_u) = F_u+ \text{Offset(2d5)}
    = F_u + \left\{
    \begin{bmatrix} -1 \\ 0 \end{bmatrix},
    \begin{bmatrix}  0 \\-1 \end{bmatrix},
    \begin{bmatrix}  0 \\ 0 \end{bmatrix},
    \begin{bmatrix}  0 \\ 1 \end{bmatrix},
    \begin{bmatrix}  1 \\ 0 \end{bmatrix}
    \right\}
\label{eq:off_2d5}
\end{equation}
Specifically, when $F_u$ locates at $C_i$'s left-up corner,
substituting $s_x=s_y=2$, the intersection of $P_\text{to}$ and $\text{Ngb}(F_u)$ contains three coordinates of possible $F_v$.
They exactly correspond to three chains in Figure~\ref{fig:Galerkin-overview}(b3), (b4) and (b5):
\begin{equation} \small
    P_{\text{to}} \cap \text{Ngb}(F_u) = C_i + \left\{
        \begin{bmatrix}  -2 \\ 1 \end{bmatrix},
        \begin{bmatrix}  -1 \\ 1 \end{bmatrix},
        \begin{bmatrix}  -1 \\ 2 \end{bmatrix}
    \right\}
\end{equation}
The other two chains in Figure~\ref{fig:Galerkin-overview}(b1) and (b2) are derived in a similar way when $F_u$ locates at $C_i$'s left side and up side, respectively.
The operator-perspective computation is equivalent to matrix-perspective computation by disassembling a nonzero value of $A^C$ into all possible pieces of contributions.

However, such derivation is too tedious to code manually. 
Table~\ref{tab:LoC} displays the numbers of influence chains for different combinations of stencil patterns of $R$, $A^F$, and $P$, including the commonly-used patterns such as 3d7, 3d19, and 3d27.
There are about 2200 influence chains when $R$, $A^F$, and $P$ are all of 3d27 pattern. 
The large numbers of chains indicate that it is too labor-intensive for humans to derive them chain-by-chain correctly and then code them line-by-line without bugs.
Therefore, the symbolic derivation could be performed by softwares or languages with symbolic-computation support (e.g., Python, Octave, Matlab, Mathematica)
, followed by C++ code generation in our work.
As listed in Table~\ref{tab:LoC}, the huge numbers of lines of codes (LoC) demonstrate the necessity of automatic analysis and generation.
The generated codes for the above 2D example in Figure~\ref{fig:Galerkin-overview} are in Listing~\ref{code} in supplemental material~\ref{txt:code} due to page limit, illustrating the complicated coding.

\begin{table}[htbp]
\vspace{-0.3em}
  \footnotesize
  \caption{Lines of codes generated in different combinations. 
  }
  \label{tab:LoC}
    \begin{threeparttable}
  \begin{tabular}{c|ccccccc}
    \toprule
    $R$  & 3d8c & 3d8c & 3d8c & 3d7v & 3d27v & 3d27v & 3d27v\\
    $A^F$& 3d7  & 3d19 & 3d27 & 3d7  & 3d15  & 3d19  & 3d27\\
    $P$  & 3d8c & 3d64c& 3d64c& 3d7v & 3d27v & 3d27v & 3d27v\\
    \midrule
    Chains& 56 & 1216 & 1728 & 37 & 1208 & 1685 & 2197 \\
    LoC & 244 & 2364 & 2916 & 168 & 2879 & 3496 & 4048 \\
  \bottomrule
\end{tabular}
\begin{tablenotes}    
\item[] Suffix 'c' in $P,R$ stands for cell-centered~\cite{Wesseling-101,Kwak_Cell_Centered} interpolation and restriction, while 'v' for vertex-centered~\cite{KhalilWesseling-103,Dendy_BlackBox}.
\end{tablenotes}
\end{threeparttable}
\vspace{-0.8em}
\end{table}

The two-step multiplications in the triple-matrix product are fused into one single kernel at the cost of redundant computation.
But the reduced memory footprint and communication volume improve the computational efficiency significantly.
The fusion strategy trends to well fit the features of current architectures.
A computation kernel could be generated for each combination of $R$, $A^F$, and $P$ of different patterns.
As shown in Listing~\ref{code}, when the matrices data are organized in SG-DIA format, the computation, memory access, and communication patterns are similar to an ordinary stencil.
The patterns inside the innermost loop are slightly more irregular (i.e., not fully contiguous) than a stencil, and therefore, software prefetching instructions are also generated and embedded at suitable positions to improve performance.
Furthermore, standard optimizations for stencils, such as tiling~\cite{stencil_35D} and vectorization~\cite{stencil_vec}, could be applied based on our code-generation tool.


We consider symbolic derivation and automatic code generation to be an indispensable key to the success of StructMG, which can extend to various architectures, such as each thread on GPU computing a coarse element.
It handles the complicated implementation for different dimensional coarsening and different combinations of $R$, $A^F$, and $P$, fulfilling the principle \textbf{P1} and \textbf{P2} at the same time.
Based on the native data structures (i.e., grids and stencils), our method avoids introducing general sparse matrix storage that is less performance-favored.
With a complete fusion of the products, it has a smaller memory footprint than two-step multiplication.

\subsection{Smoothers: Vectorized Level-Based SpTRSV}
Most often used smoothers in multigrid include Jacobi, Gauss-Seidel (GS), and ILU.
Since Jacobi is free of read-after-write dependence and could be realized by SpMV, its parallel optimization are straightforward.
Symmetric GS (including forward and backward passes) and ILU are more sophisticated. Their formulas are as follows.
\vspace{-0.em}
\begin{align*}
\text{GS-Forward}&: & u^{(t+\frac{1}{2})}=&u^{(t)}+w(D+wL)^{-1}(b-Au^{(t)}) \\
\text{GS-Backward}&: & u^{(t+1)}=&u^{(t+\frac{1}{2})}+w(D+wU)^{-1}(b-Au^{(t+\frac{1}{2})})\\
\text{ILU}&: & u^{(t+1)}=&u^{(t)}+\widetilde{U}^{-1}\widetilde{L}^{-1}(b-Au^{(t)})
\end{align*}
where $L, D, U$ are the lower triangular, diagonal, and upper triangular part of the matrix $A^{*h}$ (as denoted in Figure 2), and $\widetilde{L},\widetilde{U}$ are factorized lower and upper triangular matrix of $A^{*h}$.
Note that GS and ILU are essentially triangular solvers, and parallelizing them efficiently requires considerable effort.
Following the principle \textbf{P3}, StructMG respects the dependency order strictly to provide more robust smoothers than RB-GS in PFMG or hybrid Jacobi/GS in BoomerAMG that exploit parallelism at the cost of slower convergence.
The efficiency is obtained without slowing down convergence of StructMG because the parallel results of multiple threads are the same as the serial one.

Since dependence-preserving GS, and ILU are all special forms of SpTRSV, they can be implemented in a unified parallel framework.
A widely-used method is level scheduling~\cite{Maxim-2011}.
A directed acyclic graph derived from the triangular matrix is divided into multiple levels, and nodes on the same level can be processed in parallel.
The previous study sparsified synchronization by analysis of static scheduling and used lower-overhead P2P rather than a barrier~\cite{Park-ICS14}.
Many high-performance implementations of SpTRSV focus on unstructured problems~\cite{Wang_PPoPP18, Dufrechou_2018, Xie_icpp2021} but cannot be used directly in StructMG because they violate the principle \textbf{P1}.
Optimizations of SpTRSV on structured grids have been reported on manycore architectures~\cite{Wang-ICPP18, Zhu-SC21}, relying on specific hardware support. 
However, most ARM and X86 multicore architectures can only synchronize threads by accessing the last-level shared cache among cores. 
Due to the high latency, we must reduce the number of synchronizations. Therefore, we proposed a novel design of SpTRSV on structured grids suitable for nowadays mainstream~SMP~processors.

\begin{figure}[h] 
  \centering
  \includegraphics[width=\linewidth]{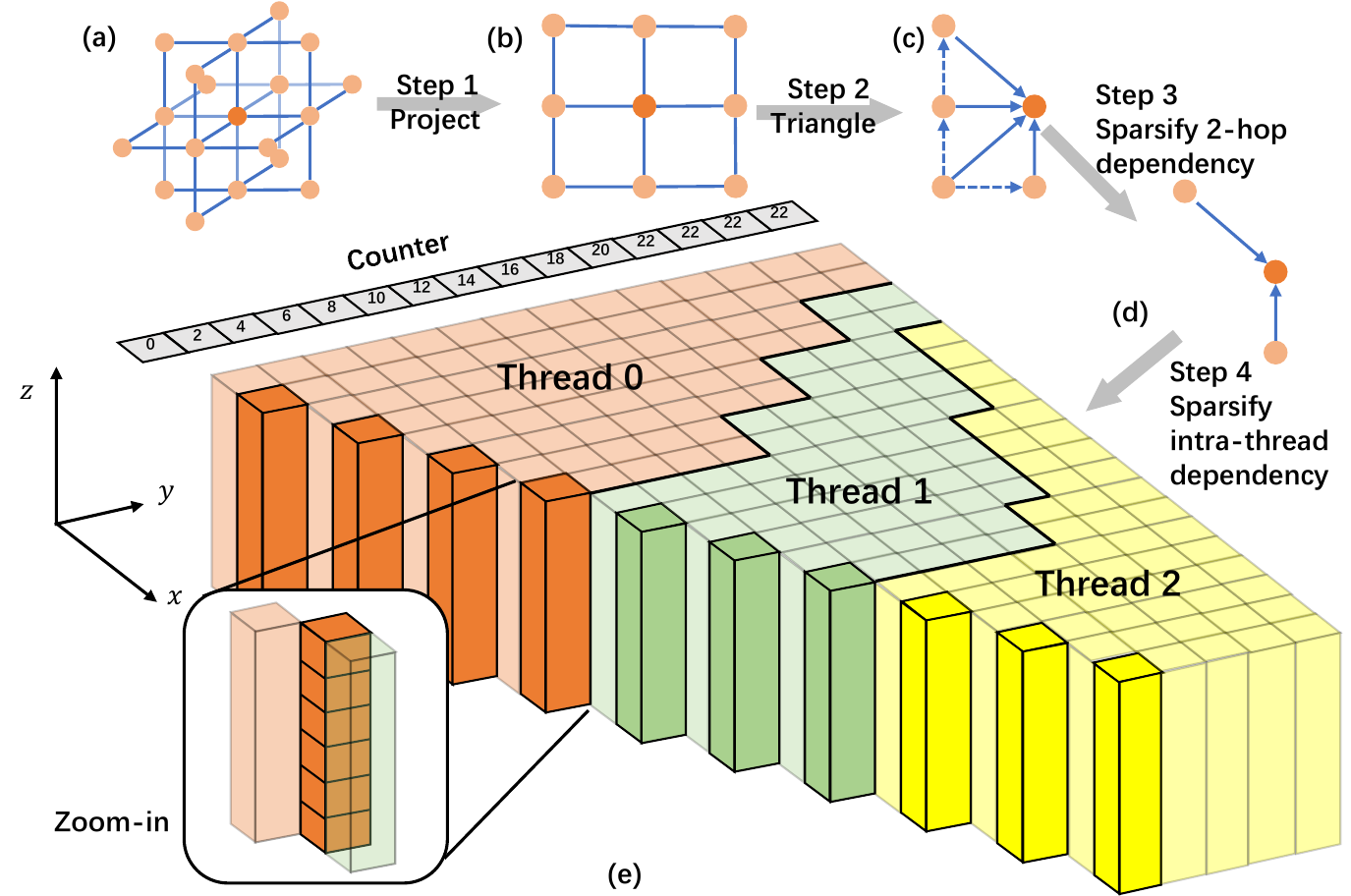}
  \caption{Illustration of SpTRSV on a structured grid.}
  \label{fig:SpTRSV}
    \vspace{-1.0em}
\end{figure}

Figure~\ref{fig:SpTRSV} illustrates how multiple threads in a process collaborate SpTRSV on structured grids.
Level-based task scheduling is applied in the outer two dimensions ($x$ and $y$) with minimal task granularity as a single column. 
A 3d* stencil is first projected into a 2d* stencil, as Step 1 shows an example of 3d19 into 2d9. 
The triangular part of 2d9 is extracted as skewed 2d5 in Figure~\ref{fig:SpTRSV}(c) after Step 2.
Ideas of sparsifying synchronization in~\cite{Park-ICS14} are employed, and such dependence analysis can be carried out statically due to the structured features of the grids.
First, 2-hop dependence is sparsified. 
Note the dependence of an element's dependence is not its direct dependence.
For example, the left and left-down neighbors' solid arrows towards the center in Figure~\ref{fig:SpTRSV}(c) are eliminated in Step 3 because the dashed arrows departing from them indicate the 2-hop dependence.
According to the cascade of dependence, they must have been finished before the computation at the center.
Second, intra-thread dependence is sparsified. Based on the static scheduling strategy, as shown in Figure~\ref{fig:SpTRSV}(e), thread 0 is responsible for the orange domain, 1 for green, and 2 for yellow. 
Columns with darker colors represent the forefront level. 
No synchronization is needed for the dependence that has been processed by the same thread on previous levels. 
For example, the orange column in the zoom-in plot only needs to wait for the green column to be finished. 
A shared counter array contains $x$ coordinates of the latest finished positions at each $y$ position.
Inter-thread synchronization by looking up a counter table is necessary only if a column adjoins the thread’s domain boundaries (marked by bold lines in Figure~\ref{fig:SpTRSV}(e)) and its corresponding dependent column lies in other threads' domains.
Counter array enables P2P synchronization instead of less-efficient barriers among all threads. 
This strategy of work distribution can~extend to other architectures, such as GPU~\cite{Struct_SpTRSV_Anonymous} with a few~adaptations.

Computation within a column is vectorized.
Elements in neighboring columns can be loaded, multiplied by matrix coefficients, accumulated, and stored to intermediate results in a vectorized format.
For PGS and ILU, where updating values is point-wise, intermediate results are extracted from the vector register, added with contributions of elements in the current column, and finally written back to memory in order, which is similar to~\cite{Zhu-SC21}.
For LGS, where values of a column are updated simultaneously, the above procedure assembles a right-hand-side vector, and a tridiagonal equation is solved within the current column.
Our tridiagonal solver adopts Thomas algorithm~\cite{TDMA} with the lowest computational complexity.

\section{Experimental Setup}

\subsection{Problems} \label{txt:problems}

Testing problems' characteristics are listed in Table~\ref{tab:problems}, including 
nonzero patterns (Base Pat. field), 
total degrees of freedom (\#dof field in a unit of million),
anisotropy strength (Aniso. field),
precision (Prec. field), 
condition numbers (Cond. field),
iterative solvers (Solver field),
converging tolerance given by users' applications (Tol. field, where 'r' stands for relative tolerance and 'a' for absolute one).
More detailed numerical distribution of the matrix and the statistics of the anisotropic metric for these problems can be found in Figure 1 and Figure 5 of~\cite{ZY_FP16}, respectively.
Directions of anisotropy are shown in Table~\ref{tab:anisotropy}.
Data are available online.\footnote{See https://zenodo.org/records/10023590 except for \cName{Laplace} because it has constant coefficients, making it straightforward for readers to construct on their own.}

\begin{table}[htbp]
  \footnotesize
  \caption{Test problems characteristics.}
  \label{tab:problems}
  \begin{threeparttable}
  \begin{tabular}
  {c|ccccccc}
    \toprule
    Problem  & Base Pat.\tnote{1}& \#dof & Aniso. & Prec. & Cond. & Solver & Tol. \\
    \midrule
    \cName{Laplace}  & 3d7  & 16.8M & None & FP64 & 2e+04 & CG   & r.<$10^{-9}$ \\
    \cName{rhd}      & 3d7  & 2.10M & Low  & FP64 & 1e+08 & CG   & r.<$10^{-9}$ \\
    \cName{weather}  & 3d19 & 637M  & High & FP32 & 1e+05\tnote{2} & GMRES& a.<$10^{-5}$  \\
    \cName{rhd-3T}   & 3d7  & 6.30M & High & FP64 & 1e+15 & CG   & r.<$10^{-9}$  \\
    \cName{oil-4C}   & 3d7  & 31.5M & High & FP64 & 1e+05 & GMRES& r.<$10^{-4}$  \\
    \cName{solid-3D} & 3d15 & 11.8M & Low  & FP64 & 1e+07 & CG   & r.<$10^{-9}$  \\
  \bottomrule
\end{tabular}
\begin{tablenotes}    
\footnotesize               
\item[1] 3d7, 3d15 and 3d19 may expand to 3d27 on coarser grids.
\item[2] A matrix of only 9.9M rows discretized from the same problem is estimated instead since the original matrix of 637M is too large to compute condition number.
\end{tablenotes}
\end{threeparttable}
\vspace{-1.2em}
\end{table}

Our selection of test problems aims to cover different domains and characteristics as extensively as possible.
Problems with a numerical suffix are vector PDE problems
, where the number indicates the \#dof (i.e., unknowns) per element.
Otherwise, they are scalar PDE problems.
%
\cName{Laplace} is an idealized benchmark problem in AMG evaluation and modeling, such as \textit{hypre}'s reports~\cite{osti_2021,FalgoutLi-130} and HPCG~\cite{Zhu-SC21}. 
It discretizes the Laplace equation with 7-point finite difference on a uniform 3D Cartesian grid and Dirichlet boundary conditions. 
It is fully isotropic and has constant coefficients.
\cName{Rhd} and \cName{rhd-3T} are from radiation hydrodynamics~\cite{XuMo2017}. 
"3T" means three temperatures (radiation, electron, and ion).
Discretized from three-temperature equations, \cName{rhd-3T} is highly-anisotropic due to non-smooth coefficients and multi-physics coupling.
\cName{Rhd} is from the decoupling of the \cName{rhd-3T} system and is weaker anisotropic.
\cName{Oil-4C} are from petroleum reservoir simulation. 
Settings of SPE1 and SPE10 benchmarks~\cite{Tenth_SPE, SPE_CSP_2} are combined to generate larger cases via OpenCAEPoro~\cite{CAEPoro_SPE10}.
"4C" means four components (oil, water, gas, and dissolved gas in live oil).
\cName{Oil-4C} is highly-anisotropic due to inhomogeneous permeability.
\cName{Weather} is from atmospheric dynamics, provided by GRAPES-MESO~\cite{GRAPES-OL, GRAPES}, the national weather forecasting system of China Meteorological Administration~\cite{CMA-OL}. 
The strong anisotropy comes from irregular earth topography and nonuniform latitudinal spacing.
The largest case of 637M \#dof is a 2km resolution of Chinese region in Dec 2018.
\cName{Solid-3D} cases are generated by ourselves.
"3D" means three displacements associated with each element.
Discretized from the weak form of linear elasticity problem in solid mechanics~\cite{Elasticity}, \cName{Solid-3D} is isotropic and has homogeneous coefficients.
\cName{Laplace} and \cName{solid-3D} are from idealized problems, and the others are from real-world applications and more complicated.

\subsection{Solvers}
Iterative solvers, CG and GMRES (restart length of 10), and multigrid preconditioners, SMG, PFMG, SysPFMG, and BoomerAMG, are all from \textit{hypre}-2.25.0~\cite{HYPRE-OL}. To be fair in comparison of different multigrids, our implementations of iterative solvers follow the codes of \textit{hypre} strictly without optimizations.
All multigrids apply V-Cycle for once with pre- and post-smoothing both for once.
Due to the limited number of parameters in SMG, PFMG, and SysPFMG, exhaustive trial-and-error search to minimize time-to-solution is feasible.
The best combinations of parameters are picked for them.
BoomerAMG is more complicated because its large number of parameters makes an exhaustive search infeasible. 
Therefore, we adopt the problem-specific tuning from domain experts and follow the settings of~\cite{XuMo2017} for radiation hydrodynamics and~\cite{CAEPoro-csrfasp} for petroleum reservoir problems.
Prior knowledge of applying unstructured AMG in our \cName{weather} cases is not available.
Based on the best practice guide summarized by the \textit{hypre} team~\cite{10.1137/130931539}, we have tried different coarsen types, interpolation types, strong thresholds, and relaxation weights to find the minimal time-to-solution, which we consider is adequate for an \textbf{auxiliary baseline}
\footnote{
We want to emphasize that StructMG is not inherently in the same track as unstructured AMG, as shown in Figure~\ref{fig:correspond_to_hypre}.
\textbf{It neither intends nor can replace unstructured AMG like BoomerAMG.}
Its goal is to provide higher-performance solutions for structured problems than the existing structured multigrids such as SMG and PFMG.
Therefore, after trying our best to find the optimal parameter combination, even though it may not be the real best in the entire parameter space, we consider it adequate for BoomerAMG as an \textbf{auxiliary baseline} for comparison.
}
.
Solid-3D is a classical elasticity problem, and BoomerAMG needs near-nullspace vectors as the additional input
\footnote{
An unknown-based approach not requiring rigid body modes could be adopted, but its \#iter increased to a great extent and resulted in a longer total time. Results refer to https://gitee.com/zongy17/baseline-mg/raw/master/AMG\_woRBM\_results.pdf
}
~\cite{elasticity_NLA} that StructMG and SysPFMG do not require.
The aggressive levels are set to 1 except for \cName{oil-4C} because of divergence.
Their detailed parameter settings are listed in Artifact Description Appendix~\ref{txt:setting}.
Particularly for \cName{Laplace}, the idealized benchmark problem in common literature, we follow the same settings of SMG, PFMG, and BoomerAMG with \textit{hypre}'s official evaluation~\cite{osti_2021}.

\subsection{Machines}
Experiments are evaluated on ARM and X86 clusters, as shown in Table~\ref{tab:machines}. 
For a fair comparison, all multigrids are tested with the same number of physical cores.
Multi-threading is also supported by \textit{hypre} and hybrid MPI/OpenMP usually performs better than pure MPI on NUMA architectures~\cite{Baker2012}. Therefore, the best result at a specific number of cores is reported out of various tests of different MPI/OpenMP ratios with load-balance process partitions.
1:1 (MPI-only), 1:2, 1:4, 1:8, 1:16, and 1:32 are tested when each NUMA has 32 available cores.
1:1, 1:2, 1:3, 1:5, 1:6, 1:10, 1:15, and 1:30 are tested when only 30 cores are available because ARM's McKernel mechanism requires a 2-core reservation for OS to reduce system noise for large-scale tests.
For example, 1:4 means each MPI process runs with 4 threads, and thus 32 processes fully utilize one 128-physical-core node.
Each node is equipped with two processors.

\begin{table}[htbp]
\vspace{-1.em}
  \footnotesize
  \caption{Machines Configurations.}
  \label{tab:machines}
  \begin{tabular}{c|cc}
    \toprule
    System & ARM & X86 \\
    \midrule
    Processor & Kunpeng 920-6426&  AMD EPYC-7H12 \\
    Frequency & 2.60 GHz & 2.60\textasciitilde 3.30 GHz \\
    Cores per node & 128 (64 per socket) & 128 (64 per socket) \\
    L1/L2/L3 per core & 64 KB/512 KB/1 MB & 32 KB/512 KB/4 MB\\
    Stream Triad BW & 138 GB/s per socket & 100 GB/s per socket\\
    Memory per Node & 512 GB DDR4-2933 & 256 GB DDR4-3200 \\
    Max Nodes/Network & 64 nodes/IB(100 Gbps) & 256 nodes/IB(100 Gbps) \\
    MPI/Compiler & OpenMPI-4.1.4/gcc-9.3.0 & OpenMPI-4.1.4/icc-19.1.3  \\
  \bottomrule
\end{tabular}
\vspace{-1.0em}
\end{table}

\section{Performance Results and Analysis}
The first three subsections of this section will respectively validate that the principles \textbf{P1} to \textbf{P3} we proposed have been meticulously fulfilled.
Strong and weak scalability tests will then be presented.

\subsection{Kernel Performance}\label{txt:abla_kernel}
The kernels in this subsection utilize all cores in a NUMA by one process (i.e., 32 threads for ARM, and 64 threads for X86). 

\subsubsection{SpMV \& SpTRSV} \label{txt:abla_kernel_spmv_trsv}
SpMV is the kernel of Jacobi smoothers, and residual computation, while SpTRSV is the kernel of most-often used smoothers such as GS and ILU. 
Particularly, a representative profiling~\cite{OptimMG_multiCores} indicated that GS smoother (a specialized form of SpTRSV) accounts for 78\% of the entire execution time of the HPCG benchmark, while SpMV has the second-largest contribution of 20\%.
Therefore, it is necessary to make these two kernels efficient for a high-performance multigrid.
As they are typical memory-bounded kernels, their effective memory bandwidth is shown in Figure~\ref{fig:kernel_all}(a-1), all achieving levels close to the STREAM benchmark, which is considered as the upper performance limit of memory access.
Furthermore, we compare the structured kernels with the SOTA inspector-executor style routines in ARM Performance Library (ARMPL, 23.10) on ARM and MKL (2023.2) on X86.
After excluding the time of symbolic analysis in these baselines, our optimized SpMV could be 1.85x and 1.14x faster than ARMPL and MKL, respectively, in geometric average.
Our optimized SpTRSV could be 3.49x and 2.23x faster than ARMPL and MKL, respectively, as shown in Figure~\ref{fig:kernel_all}(a-2).

\subsubsection{Triple-matrix product}
\begin{figure}[h] 
  \centering
  \includegraphics[width=\linewidth]{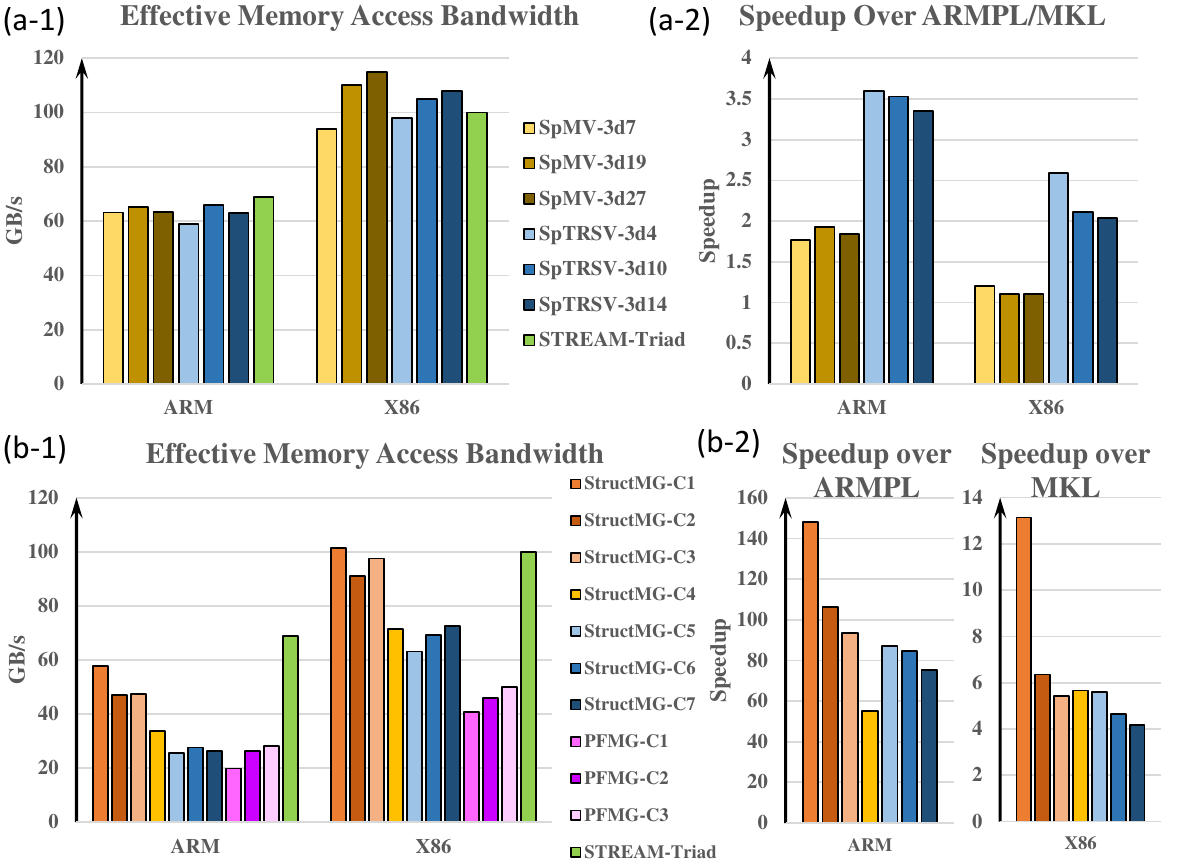}
  \caption{(a) SpMV \& SpTRSV and (b) triple-matrix product performance on two platforms.}
  \Description{...}
  \label{fig:kernel_all}
  \vspace{-1.2em}
\end{figure}

This is the critical kernel in the setup phase.
Besides being significantly faster than directly invoking the SpGEMM in ARMPL and MKL for twice, as shown in Figure~\ref{fig:kernel_all}(b-2), we also compare the effective memory bandwidth between our generated code and the manually written code in PFMG in (b-1).
The legends ending with '-C' suffix of StructMG correspond to the seven combinations of $R$, $A^F$, $P$ in Table~\ref{tab:LoC}.
In PFMG, only three combinations are available for 3D problems, all utilizing 1D coarsening with much fewer lines of codes than ours.
Despite the substantially increased code complexity due to multi-dimensional coarsening, our generated code consistently outperforms the manually written code in PFMG on both platforms.

\subsection{Complexity Analysis}\label{txt:complexity_analysis}

The grid complexity $C_{\text{G}}$ and operator complexity $C_{\text{O}}$ of different multigrids in different problems are summarized in Table~\ref{tab:complexity}, where PFMG (used for scalar PDE) and SysPFMG (used for vector PDE) share the same column in the table.
StructMG has the lowest $C_{\text{G}}$ and $C_{\text{O}}$ among structured multigrids across all problems, with a notable reduction of \textbf{69\%} in $C_{\text{O}}$ compared to SysPFMG.
Multi-dimensional coarsening helps StructMG reduce $C_{\text{G}}$ as low as aggressive coarsening in unstructured AMG.
It is worth mentioning that the non-Galerkin coarsening of PFMG can maintain the patterns of coarser-level matrices as 3d7, but it converges too slowly in our real-world problems and is thus not enabled.
As Principle \textbf{P2} has stated in Section~\ref{txt:principles}, lower complexity is favored by both speed and scalability, as it reduces not only computational workload and communication volume, but also decreases the number of communications through fewer AMG levels.
As a comparison, SMG and (Sys)PFMG require 3 to 4 times as many levels as StructMG and BoomerAMG.

\begin{table}[htbp] \footnotesize
\caption{Complexity of different multigrids.}
  \label{tab:complexity}


\begin{tabular}{cc|cccc}
     \toprule
     & $C_\text{G} / C_\text{O} $ & StructMG & SMG & (Sys)PFMG & BoomerAMG \\
     \midrule
     \multirow{3}*{\makecell{Scalar\\PDE}} 
     & \cName{Laplace}    & 1.14 / 1.55 & 1.98 / 3.11 & 2.00 / 4.29 & 1.08 / 1.33 \\
     & \cName{rhd}        & 1.14 / 1.41 & 1.98 / 3.11 & 2.00 / 4.29 & 1.08 / 1.33 \\
     & \cName{weather}    & 1.31 / 1.44 & 2.00 / 2.42 & 2.00 / 2.42 & 1.28 / 1.62 \\
     \hline
     \multirow{3}*{\makecell{Vector\\PDE}} 
     & \cName{rhd-3T}     & 1.14 / 1.55 & \multirow{3}*{\makecell{Not\\Usable}}& 2.00 / 6.04 & 1.04 / 1.14 \\
     & \cName{oil-4C}     & 1.14 / 1.14 &           & 2.00 / 3.57 & 1.87 / 2.62 \\
     & \cName{solid-3D}   & 1.14 / 1.26 &           & 2.00 / 3.60 & 1.14 / 1.75 \\
     \bottomrule
\end{tabular}


\vspace{-1.0em}
\end{table}

\subsection{Convergence Comparison}\label{txt:abla_convergence}
After proving the efficiency of smoothers in Section~\ref{txt:abla_kernel_spmv_trsv}, this subsection further examines their convergence properties.
When the multigrid configuration and process partitioning are fixed, the \#iter needed to converge using different smoothers are displayed in Table~\ref{tab:anisotropy}.
Stricter tolerances than Table~\ref{tab:problems} are adopted for more evident comparison in this subsection.
Point GS is already good enough for \cName{Laplace}, an isotropic problem, while more sophisticated smoothers show advantages for the other problems. 
When the anisotropic direction is grid-aligned, line GS is better than ILU.
When the anisotropic direction is irregular or singular, ILU significantly outperforms the other smoothers.
Particularly, as the anisotropic strength increases, the convergence gap between point GS and the best smoother widens.
These results align with Principle \textbf{P3} that structured AMG employing fix-stride coarsening requires more robust smoothers to manage the multigrid seesaw effectively, especially for anisotropic problems.
To further verify that the best smoother itself does not ensure good convergence, we also report the \#iter when using the best smoother directly as a preconditioner (i.e., without MG), which far exceed the optimal \#iter.

\begin{table}[htbp]  \footnotesize
  \caption{Ablation study of smoother effect on anisotropy. The \#iter to converge are shown (fewer is better). 
  }
  \label{tab:anisotropy}
  \begin{threeparttable}
  \begin{tabular}
  {c|cc|ccc|c}
    \toprule
    \multirow{2}*{Problem} & \multicolumn{2}{c|}{Anisotropy} & \multicolumn{3}{c|}{Different Smoothers\tnote{1}} & Best Smoother \\
      & Strength & Direction & PGS & LGS & ILU & \textbf{Without} MG \\
    \midrule
    \cName{Laplace}  & None & None           & \underline{\textbf{13}} & 13 & 13 & 287 \\
    \cName{rhd}      & Low  & Irregular      & 45 & 45 & \underline{\textbf{22}} & 511  \\
    \cName{solid-3D} & Low  & Aligned/Skewed\tnote{2} & 15 & 14 & \underline{\textbf{12}} & 150 \\
    \cName{rhd-3T}   & High & Irregular      & 90 & 82 & \underline{\textbf{25}} & 992 \\
    \cName{oil-4C}   & High & Aligned        &1200& \underline{\textbf{31}} & 43 & 295 \\
    \cName{weather}  & High& Aligned        & 66 &  \underline{\textbf{7}} &  8 &  24 \\
    \cName{weather} (g)\tnote{3}& High& Singular\tnote{4}& 60 & 21 & \underline{\textbf{13}} & 60 \\
  \bottomrule
\end{tabular}
\begin{tablenotes}    
\footnotesize               
\item[1] Different smoothers with the same StructMG configuration to observe smoother effect.
\item[2] Aligned in the first three levels and skewed in the following coarser levels.
\item[3] This is a global version of the \cName{weather} problem. Although it is actually solved by Semi-StructMG~\cite{Semi-StructMG},
it is listed here to facilitate comparative analysis with the \cName{weather} case.
\item[4] Direction changes complicatedly near the polar regions since the longitude grid lines cluster at the North and South Poles, exhibiting non-Cartesian singular characteristics.
\end{tablenotes}
\end{threeparttable}
\vspace{-1.0em}
\end{table}

\subsection{Strong Scalability}
\begin{figure*}[htbp] 
  \centering
  \includegraphics[width=\textwidth]{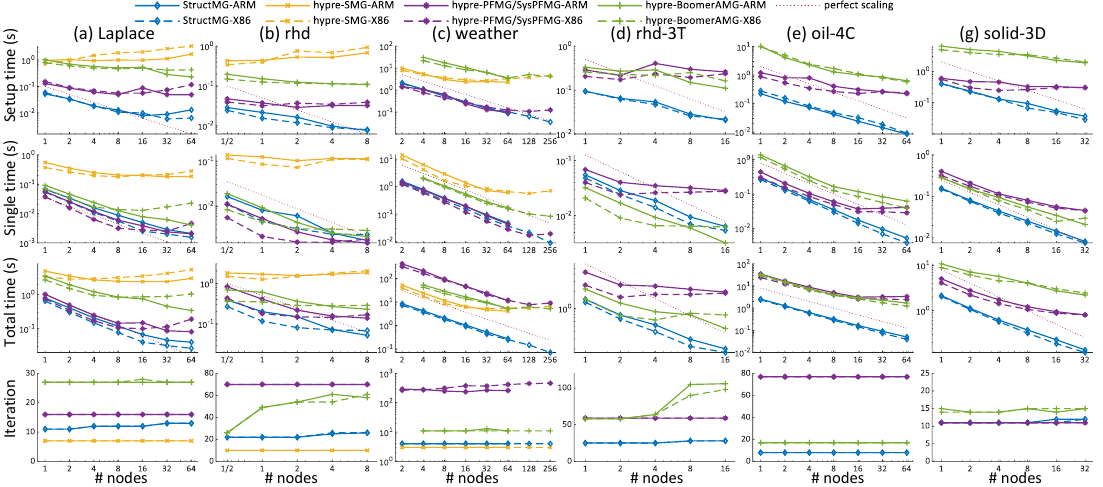}
  \caption{Strong scalability results. Lower is better in all subfigures. Subfigures share legends at the top, y-axis labels at the left, and x-axis labels at the bottom. Four rows of subfigures show $T_\text{setup}$, $T_\text{single}$, $T_\text{tot}$, and \#iter (denoted in Equation~(\ref{eq:time})), respectively.}
  \label{fig:strong}
  \vspace{-1.0em}
\end{figure*}

The preceding subsections have shown that StructMG maximizes the performance benefits for structured problems in multiple aspects, pushing its potential to the limit.
The end-to-end tests will further verify its advantage.
Strong scalability tests are reported in Figure~\ref{fig:strong}.
The breakdown of the terms in Equation~(\ref{eq:time}) is shown in four subplots arranged in rows.
PFMG (used for scalar PDE) and SysPFMG (used for vector PDE) share the same legend.

The results on ARM are summarized first.
StructMG exhibits the best scalability in both setup and solve phases for both types of PDE problems.
In \cName{rhd}, the smallest-size scalar PDE problem that best reflects scalability, StructMG scales to the maximal degree of parallelism with an efficiency of \textbf{46\%} in terms of $T_\text{tot}$. In contrast, SMG, PFMG, and BoomerAMG could only scale with 6.2\%, 36\%, and 19\%, respectively.
Similarly, in \cName{rhd-3T}, the smallest-size vector PDE problem, StructMG, SysPFMG, and BoomerAMG scale to the maximal degree of parallelism with efficiencies of \textbf{45\%}, 14\%, and 15\%, respectively.
SMG is the least scalable multigrid that always reaches the saturated point of performance earliest.
Even if it has the fewest \#iter, SMG is the slowest regardless of degrees of parallelism.
Its too long $T_\text{single}$ makes it competitive in $T_\text{tot}$ only when anisotropy is so strong that the others have excessive \#iter, as in \cName{weather}.
PFMG scales as perfectly as StructMG in \cName{weather} with the largest-size \#dof, but degrades in the other problems.
PFMG always has more \#iter than StructMG because its smoothers (weighted Jacobi or RB-GS) are less convergence-favored than dependence-preserving GS and ILU.
Such difference of \#iter is huge in \cName{weather}, as plotted logarithmically in (c-4).
SysPFMG has similar \#iter with StructMG but is less scalable and much slower.
The gaps between the scaling curves of SysPFMG and StructMG widen quickly with increasing parallelism in these problems.
BoomerAMG has similar scaling behaviors with PFMG and SysPFMG.
But noticeably, it suffers increases in \#iter when the number of cores increases in the smallest-size problems such as \cName{rhd} and \cName{rhd-3T}.
This is because when the problem size is fixed, the effect of hybrid Jacobi/GS smoother gradually degenerates into Jacobi-only as the number of cores increases.
It is worth mentioning that unstructured AMG's too-heavy setup makes it sometimes slower than SysPFMG in $T_\text{tot}$, even if it has faster $T_\text{single}$ and similar \#iter.

The scaling results on X86 are similar to ARM.
With more available X86 nodes, the scalability advantage of StructMG is further revealed in \cName{weather}.
StructMG continues scaling to 256 nodes with an efficiency of \textbf{80\%}, but SMG, BoomerAMG, and PFMG have shown saturated performance since 16, 64, and 128 nodes, respectively.

\subsection{Weak Scalability}
\begin{figure}[htbp] 
  \centering
  \includegraphics[width=\linewidth]{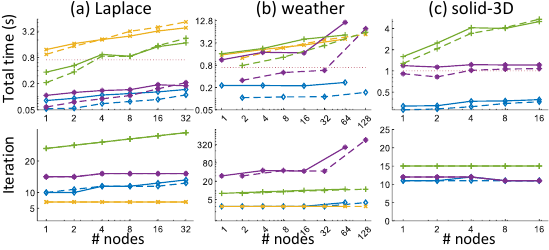}
  \caption{Weak scalability results. Legends match Figure~\ref{fig:strong}.}
  \label{fig:weak}
  \vspace{-1.0em}
\end{figure}

Three problems are evaluated for weak scalability and shown in Figure~\ref{fig:weak}.
Only total time and \#iter are shown due to page limit.

For \cName{Laplace} cases that are commonly used as benchmarks, StructMG, SMG, PFMG, and BoomerAMG scale to 32 nodes with efficiencies of \textbf{48\%}, 18\%, 28\%, and 9.6\% on X86, respectively.
Our testing results of \textit{hypre}'s multigrids for \cName{Laplace} are similar to those in previous literature~\cite{osti_2021, Baker2012}.
For \cName{weather}, StructMG scales to 128 nodes on X86 with an efficiency of \textbf{72\%}, which drops from \textbf{94\%} at 32 nodes due to an increase in \#iter from 3 to 4. 
Its weak scaling behavior of single-iteration time is near-perfect. 
In contrast, SMG scales with 23\% even if it maintains a constant \#iter of 3, 
while BoomerAMG scales with 13\% at the cost of moderately increasing \#iter. 
PFMG degrades slightly in single-iteration solve time but significantly in \#iter, resulting in a parallel efficiency of only 4\% at 128 nodes.
SMG, PFMG, and our StructMG all employ fix-stride coarsening, but PFMG has significantly more \#iter, which implies the weakness of point-wise weighted Jacobi when solving larger-size anisotropic problems.
Noticeably, by applying more sophisticated LGS, StructMG could reach a convergence rate close to SMG. 
For \cName{solid-3D} from vector PDE, StructMG, SysPFMG, and BoomerAMG scale to the maximal degree of parallelism with efficiencies of \textbf{77\%}, 72\%, and 23\%, respectively.
Results on ARM are similar to those on X86.

\subsection{Time-to-Solution}
The minimal, geometric average, and maximal time-to-solution speedups of StructMG over others with the same degree of parallelism in all test runs are summarized in Table~\ref{tab:speedup}.
Owing to StructMG's scalability advantage, the minimal speedup usually occurs at the minimal degree of parallelism, while the maximal speedup occurs at the maximal degree of parallelism.
StructMG obtains average speedups of \textbf{15.5x} over SMG, \textbf{5.5x} over PFMG, \textbf{6.7x} over SysPFMG, and \textbf{7.3x} over BoomerAMG in total time across all test runs.
The prominent speedups show that StructMG has potential to work as an alternative to SMG, PFMG, and SysPFMG.

\begin{table}[htbp] \footnotesize
\vspace{-0.6em}
  \caption{Min/Avg/Max time-to-solution speedups of StructMG over other multigrids in all test runs.}
  \label{tab:speedup}
\begin{tabular}{c|ccccccc}
    \toprule
    Min / Avg / Max & & SMG & (Sys)PFMG & BoomerAMG \\
    \midrule
\multirow{2}{*}{\cName{Laplace}}& ARM & 6.18 / 21.4 / 74.8& 1.20 / 1.53 / 2.23 & 4.44 / 7.02 / 11.2\\
                        & X86 & 5.15 / 36.8 / 204 & 1.14 / 2.01 / 7.23 & 4.14 / 11.9 / 39.4\\
\hline
\multirow{2}{*}{\cName{rhd}}    & ARM &4.61 / 12.5 / 34.2 & 1.45 / 2.04 / 2.65 & 1.77 / 2.90 / 4.37\\
                        & X86 &5.54 / 15.1 / 29.2 & 1.55 / 1.87 / 2.44 & 1.37 / 2.97 / 4.15\\
\hline
\multirow{2}{*}{\cName{weather}}& ARM & 5.64 / 7.39 / 14.9 & 41.8 / 45.9 / 50.0 & 12.7 / 15.8 / 19.1\\
                        & X86 & 4.73 / 12.7 / 92.1 & 44.2 / 53.6 / 127  & 11.6 / 22.5 / 71.9\\
\hline
\multirow{2}{*}{\cName{rhd-3T}} & ARM & \multirow{6}*{\makecell{SMG Not\\Usable for\\Vector PDE\\Problems}} &  2.94 / 5.05 / 9.68 & 1.50 / 1.90 / 2.66 \\
                        & X86 &           &1.97 / 4.47 / 10.7 & 1.14 / 2.03 / 4.47\\
\cline{1-2}\cline{4-5}
\multirow{2}{*}{\cName{oil-4C}} & ARM &           &12.8 / 20.9 / 67.1 & 11.9 /16.8 / 33.8 \\
                        & X86 &           &10.1 / 18.4 / 61.7 & 11.4 /16.3 / 31.8 \\
\cline{1-2}\cline{4-5}
\multirow{2}{*}{\cName{solid-3D}}&ARM &           &2.35 / 3.30 / 5.92 & 5.01 / 9.19 / 16.3 \\
                        & X86 &           &1.99 / 3.21 / 6.73 & 4.31 / 9.72 / 20.5 \\
     \bottomrule
\end{tabular}
\vspace{-1.5em}
\end{table}

\section{Conclusion}
In this paper, we present StructMG, a fast and scalable AMG for structured grids.
The innovative optimization is to leverage the characteristics of structured grids to enhance kernel performance and design multi-dimensional coarsening to reduce complexity, while pairing them with robust-yet-efficient smoothers to handle complicated and anisotropic problems.
These findings can also benefit other structured AMG.
To demonstrate StructMG's effectiveness, we tested six representative problems on ARM and X86 platforms.
StructMG outperformed \textit{hypre}'s SMG, PFMG, SysPFMG, and BoomerAMG in both scalability and time-to-solution, achieving average speedups of \textbf{15.5x}, \textbf{5.5x}, \textbf{6.7x}, and \textbf{7.3x}, respectively.
StructMG also significantly improved parallel efficiencies in all problems.






\bibliographystyle{ACM-Reference-Format}
\bibliography{references}

\clearpage
\appendix

\section{Artifact Description}
\label{txt:AD}
\textbf{Persistent ID}: https://zenodo.org/records/10346358

The artifacts include code repositories of StructMG of scalar PDE version, StructMG of vector PDE version, and baseline programs which invoke \textit{hypre}'s multigrids including SMG, PFMG, SysPFMG and BoomerAMG.
These repositories are archived in zipped files in the above persistent ID.

Instructions on how to compile and run are included in these repositories.
A StructMG example of 3d27 Laplacian equation is also shown, which is the same problem of classical AMG benchmark (https://github.com/LLNL/AMG). 

\subsection*{Artifact: StructMG}
\subsubsection*{StructMG of scalar PDE version}

This artifact is StructMG's implementation for scalar PDEs. 
It solves the three scalar PDE problems in the article including Laplace, rhd, and weather.
Different branches in the repository correspond to different problems. The zipped file \texttt{smg-master.zip} is for Laplace, and weather cases, while \texttt{smg-laser.zip} is for rhd case.
They are different branches of the original repository.
Programs on ARM and X86 platform share the same codes files, and differentiate by defining macro \texttt{\_\_aarch64\_\_} for those ARM-specific primitives.

\subsubsection*{StructMG of vector PDE version}

This artifact is StructMG's implementation for vector PDEs. 
It solves the three vector PDE problems in the article including rhd-3T, oil-4C, and solid-3D.
There are two branches in the repository. 
The zipped file \texttt{sys-smg-master.zip}, corresponding to branch \texttt{master}, is for rhd-3T and solid-3D, while \texttt{sys-smg-oil-4c.zip} is for oil-4C.
Other specifications are the same as the scalar PDE version.

\subsection*{Artifact: Baselines to invoke \textit{hypre}'s multigrids}
In addition to the above persistent ID, this artifact can also be found at https://gitee.com/zongy17/baseline-mg.

This artifact includes source code files to invoke different kinds of \textit{hypre}'s multigrids as baselines in our comparisons.
File \texttt{base\_struct.c} invokes SMG and PFMG through Struct interface.
File \texttt{base\_sstruct.c} invokes SysPFMG through SStruct interface.
File \texttt{base\_IJ.c} invokes BoomerAMG through IJ interface.
The specific parameters to invoke different multigrids for different problems could be found in these files with calls of \texttt{HYPRE\_*} functions.

\subsection*{Artifact Settings of Different Multigrids} \label{txt:setting}
For Laplace, the idealized benchmark problem in common literature, we follow the same settings of SMG, PFMG, and BoomerAMG with \textit{hypre}'s official evaluation~\cite{osti_2021}.
The parameter settings of different multigrids for the six problems are listed in Table~\ref{tab:mg_settings}, as a supplement of Table~\ref{tab:problems} in the text.

\begin{table*}[htbp]
  \footnotesize
  \caption{Test problems characteristics and settings of different multigrids.}
  \label{tab:mg_settings}
  \begin{threeparttable}
  \begin{tabular}
  {c|c|c|c|c|c|c}
    \toprule
      \textbf{Case} & \textbf{Laplace} & \textbf{rhd} & \textbf{weather} & \textbf{rhd-3T} & \textbf{oil-4C} & \textbf{solid-3D} \\
    \midrule
    Base Pat.\tnote{1} & 3d7 & 3d7 & 3d19 & 3d7 & 3d7 & 3d15 \\
    \#dof& 2.10M & 31.5M & 637M & 6.30M & 31.5M & 11.8M \\
    Prec.& FP64 & FP64 & FP32 & FP64 & FP64 & FP64 \\
    Solver & CG & CG & GMRES & CG & GMRES & CG \\
    Tol. & $\Vert r\Vert_2/\Vert b\Vert_2<10^{-9}$ & $\Vert r\Vert_2/\Vert b\Vert_2<10^{-9}$ & $\Vert r\Vert_2<10^{-5}$ & $\Vert r\Vert_2/\Vert b\Vert_2<10^{-9}$ & $\Vert r\Vert_2/\Vert b\Vert_2<10^{-4}$ & $\Vert r\Vert_2/\Vert b\Vert_2<10^{-9}$ \\
    \hline
    \makecell[l]{StructMG\\settings}& 
    \makecell[l]{Enable ZeroGuess.\\3D full coarsening.\\Sym PGS smoother.\\ Cell\_3d8 interp.\\Cell\_3d8 restrict.\\ RC3d8\_A*\_PC3d64\\Relax weights: 1.0} &
    \makecell[l]{Enable ZeroGuess.\\3D full coarsening.\\ILU(0) smoother.\\ Cell\_3d8 interp.\\Cell\_3d8 restrict.\\ RC3d8\_A*\_PC3d8SIp\\(smoothed inter-\\polation weights: 0.7)} &
    \makecell[l]{Enable ZeroGuess.\\2D semi-coarsening.\\Sym LGS smoother.\\ Cell\_2d16 interp\\Cell\_2d16 restrict.\\RC2d4\_A*\_PC2d16\\Relax weights: 1.0} &
    \makecell[l]{Enable ZeroGuess.\\3D full coarsening.\\ILU(0) smoother.\\ Cell\_3d64 interp.\\Cell\_3d8 restrict.\\RC3d8\_A*\_PC3d64} &
    \makecell[l]{Enable ZeroGuess.\\3D full coarsening.\\Sym LGS  smoother.\\ Cell\_3d8 interp.\\Cell\_3d8 restrict.\\RC3d8\_A*\_PC3d8\\Relax weights: 1.0} &
    \makecell[l]{Enable ZeroGuess.\\3D full coarsening.\\ILU(0) smoother.\\ Cell\_3d64 interp\\Cell\_3d8  restrict.\\ RC3d8\_A*\_PC3d64} \\
    \hline
    
    \makecell[l]{SMG\\settings} &
    \makecell[l]{Enable ZeroGuess.\\1D semi-coarsening.} &
    \makecell[l]{Enable ZeroGuess.\\1D semi-coarsening.} &
    \makecell[l]{Enable ZeroGuess.\\1D semi-coarsening.} &
    \multicolumn{3}{c}{Not usable for vector PDE problems.} \\
    
    \hline
    \makecell[l]{PFMG /\\SysPFMG\\settings} &
    \makecell[l]{Enable ZeroGuess.\\1D semi-coarsening.\\RAP type: 0 (Galerkin).\\Relax type: 1 (weighted\\Jacobi).\\Others default.} & 
    \makecell[l]{Enable ZeroGuess.\\1D semi-coarsening.\\RAP type: 0 (Galerkin).\\Relax type: 1 (weighted\\Jacobi).\\Others default.} & 
    \makecell[l]{Enable ZeroGuess.\\1D semi-coarsening.\\RAP type: 0 (Galerkin).\\Relax type: 1 (weighted\\Jacobi).\\Others default.} & 
    \makecell[l]{Enable ZeroGuess.\\1D semi-coarsening.\\Relax type: 2 (Red-\\Black GS)\\Others default.} &
    \makecell[l]{Enable ZeroGuess.\\1D semi-coarsening.\\Relax type: 0 (Jacobi) \\Others default.} & 
    \makecell[l]{Enable ZeroGuess.\\1D semi-coarsening.\\Relax type 2 (Red-\\Black GS)\\Others default.} \\
    
    \hline
    \makecell[l]{Boomer-\\AMG\\settings} &
    \makecell[l]{Coarsen type: 8 (PMIS).\\Strong threshold: 0.25\\Interp type: 6 (ext+i)\\Max rowsum: 0.8\\ Aggressive levels: 1\\Relax type: 6 (Sym\\G.S./Jacobi hybrid)\\Others default.} &
    \makecell[l]{Coarsen type: 8 (PMIS).\\Strong threshold: 0.25\\Interp type: 6 (ext+i)\\Max rowsum: 0.9\\Trunc. factor: 0.1\\ Aggressive levels: 1\\Relax type: 6 (Sym\\G.S./Jacobi hybrid)\\Others default.} &
    \makecell[l]{Coarsen: 10 (HMIS)\\Strong threshold: 0.25\\Interp type: 6 (ext+i)\\Pmax elemts: 5\\Aggressive levels: 1\\Relax type: 6 (Sym\\G.S./Jacobi hybrid).\\Others default} &
    \makecell[l]{Coarsen: 8 (PMIS)\\Strong threshold: 0.25\\Interp type: 6 (ext+i)\\Max rowsum: 0.9\\Aggressive levels: 1\\Relax type: 6 (Sym\\G.S./Jacobi hybrid).\\Others default.} &
    \makecell[l]{NumFunctions: 4\\SetNodal: 2\\Coarsen: 6 (Falgout)\\Relax type: 26 (Sym\\G.S./Jacobi hybrid)\\Relax order: 1 (C-F)\\Others default.} &
    \makecell[l]{NumFunctions: 3\\SetNodal: 4\\Nodal Diag: 1\\ Coarsen: 8 (PMIS)\\Aggressive levels: 1\\Interp type: 6 (ext+i)\\Pmax: 6. Relax\\type: 16 (Chebyshev)\\InterpVecVariant: 2\\InterpVecQMax 4\\Strong threshold 0.25\\Others default.} \\
  \bottomrule
\end{tabular}
\begin{tablenotes}    
\footnotesize               
\item[1] 3d7, 3d15 and 3d19 expand to 3d27 on coarser grids in StructMG.
\end{tablenotes}
\end{threeparttable}
\end{table*}

\section{Generated Codes of the Example in Figure~\ref{fig:Galerkin-overview}}
\label{txt:code}
The generated codes for the 2D example in Figure~\ref{fig:Galerkin-overview} are in Listing~\ref{code}, in which stencil numberings of 2d5 and 2d9 follow the orderings of Offsets(2d5) in equation~(\ref{eq:off_2d5}) and Offsets(2d9) in equation~(\ref{eq:off_2d9}), respectively. 
To be generalizable, operator-dependent interpolation~\cite{Dendy_BlackBox} is applied where coefficients of $P$ and $R$ may vary from element to element.
The comments in Listing~\ref{code} relate the codes with the situations in Figure~\ref{fig:Galerkin-overview}(b), which helps readers to comprehend.
As shown in Listing~\ref{code}, 
the stencil-based triple-matrix product's computation, memory access, and communication patterns are similar to an ordinary stencil. 
Standard optimizations for stencils, such as tiling~\cite{stencil_35D} and vectorization~\cite{stencil_vec}, could be applied based on our code-generation tool.

\begin{lstlisting}[caption={Generated code (excerpt) of the 2D example in Figure~\ref{fig:Galerkin-overview}.},
label={code},
escapechar=\%,
basicstyle=\scriptsize ,
mathescape,
morekeywords={idx_t, data_t}]
for (idx_t X=Xbeg; X<Xend; X++)
for (idx_t Y=Ybeg; Y<Yend; Y++) {
  data_t* Rv = RC[C_IDX(X,Y)];// restrict matrix
  ... // other coarse neighbors
  if (CHECK_BDR(X-1,Y+1)) {
    data_t res = 0.0;// accumulated result
    data_t* Pv = PC[C_IDX(X-1,Y+1)];// prolong matrix
    data_t* ptr;// pointer to fine-grid matrix
    {// when u=(2*X+bx-hx-1,2*Y+by-hy)
      data_t tmp = 0.0;
      ptr = AF[F_IDX(2*X+bx-hx-1,2*Y+by-hy)];
      tmp += Pv[6] * ptr[3];// Figure (b1)
      res += tmp * Rv[1];
    }
    {// when u=(2*X+bx-hx-1,2*Y+by-hy+1)
      data_t tmp = 0.0;
      ptr = AF[F_IDX(2*X+bx-hx-1,2*Y+by-hy+1)];
      tmp += Pv[6] * ptr[2];// Figure (b3)
      tmp += Pv[7] * ptr[3];// Figure (b4)
      tmp += Pv[3] * ptr[0];// Figure (b5)
      res += tmp * Rv[2];
    }
    {// when u=(2*X+bx-hx,2*Y+by-hy+1)
      data_t tmp = 0.0;
      ptr = AF[F_IDX(2*X+bx-hx,2*Y+by-hy+1)];
      tmp += Pv[6] * ptr[0];// Figure (b2)
      res += tmp * Rv[5];
    }
    // write to coarse-grid matrix
    AC[C_IDX(X,Y)][2] = res;// neighbor index: 2
  }
  ...// other coarse neighbors
}
\end{lstlisting}

\end{document}